\documentclass[a4paper,12pt]{article}
\usepackage{amsmath,amssymb}
\usepackage{color}
\usepackage{epsfig}
\usepackage{graphicx}

\definecolor{green1}{rgb}{0.1,0.65,0}
\definecolor{green2}{rgb}{0.1,0.6,0.1}
\definecolor{blue1}{rgb}{0.14,0.6,1.0}
\definecolor{blue3}{rgb}{0.05,0.05,0.5}
\definecolor{viol}{rgb}{0.4,0,0.9}
\definecolor{black}{rgb}{0,0,0}

\definecolor{cola}{rgb}{0.4,0,0.4}
\definecolor{colb}{rgb}{0,0.3,0.4}
\definecolor{colc}{rgb}{0.4,0.3,0}
\definecolor{oran}{rgb}{.95,.35,0}
\definecolor{brass}{rgb}{0.6,0.15,0}
\definecolor{whit}{rgb}{1,1,1}
\definecolor{grey1}{rgb}{.6,.6,.6}

\def\bla#1{\textcolor{black}{#1}}

\newtheorem{theorem}{Theorem}  [section]
\newtheorem{proposition}[theorem]{Proposition}
\newtheorem{lemma}[theorem]{Lemma}

\newtheorem{remark}[theorem]{Remark}

\numberwithin{equation}{section}

\font \lgr cmmib10 scaled \magstep1

\newfont{\bcal}{cmbsy10 scaled \magstep1}
\newfont{\ccal}{cmsy10 scaled \magstep1}
\newfont{\ctv}{msam10}

\newcommand{\bbox}{\mbox{\ctv \symbol{4}}}
\def\QED{{$\hfill\bbox$}}
\newenvironment{pf}[1]{\par\vskip1mm{\noindent\it #1.}\ }{\QED\par\vskip2mm}

\newcommand{\real}{\mathbb{R}}

\newcommand{\tens}{\mathbb{T}_{\rm sym}^{3\times 3}}

\def\dive{\mbox{\rm div\,}}
\def\sign{\mbox{\rm sign}}

\def\supess{\mathop{\mbox{sup\,ess}\,}}

\def\expe{{\rm e}}

\def\mae{\mbox{a.e.}}
\def\pif{\mbox{ if }}

\def\vect#1#2{\left(\!\begin{array}{ll} #1\\#2\end{array}\right)\!}

\def\bfsi{\mbox{\lgr \char27}}
\def\bfe{\mbox{\lgr \char34}}

\def\bfde{\mbox{\lgr \char14}}

\def\tbfe{\tilde\bfe}

\newcommand{\io}{\int_\Omega}
\newcommand{\ipo}{\int_{\partial\Omega}}

\def\bfA{\mathbf{A}}

\def\bff{\mathbf{f}}

\def\bfn{\mathbf{n}}
\def\bfq{\mathbf{q}}
\def\bfu{\mathbf{u}}

\def\bfdd{\mbox{\bf \,:\,}}

\def\dd{\mathrm{\,d}}

\def\barr{\begin{array}}
\def\earr{\end{array}}

\def\bpf{\begin{pf}}
\def\epf{\end{pf}}

\numberwithin{equation}{section}

\begin{document}

\title{\bf Phase separation in a gravity field
}
\date{}

\author{
\renewcommand{\thefootnote}{\!$\fnsymbol{footnote}$}
Pavel Krej\v{c}\'{\i}
\footnote{Mathematical Institute,
 Academy of Sciences of the Czech Republic, \v{Z}itn\'a 25, CZ-11567 Praha 1, Czech Republic,
E-mail {\tt  krejci@math.cas.cz}},
, Elisabetta Rocca
\footnote{Dipartimento di Matematica, Universit\`a di Milano,
Via Saldini 50, 20133 Milano, Italy, E-mail {\tt elisabetta.rocca@unimi.it}}
, and J\"urgen Sprekels
\footnote{Weierstrass Institute for Applied
Analysis and Stochastics, Mohrenstr.~39, D-10117 Berlin,
Germany, E-mail {\tt  sprekels@wias-berlin.de}}
}

\maketitle
\date{}


\vspace{-.4cm}

\noindent {\bf Abstract.} We prove here well-posedness and convergence to equilibria
for the solution trajectories associated to a model for solidification of a liquid
content of a rigid container in a gravity field.
We observe that the gravity effects, which can be neglected without considerable changes of the
process on finite time intervals, have a
substantial influence on the long time behavior of the evolution
system. Without gravity, we find a temperature interval, in
which all phase distributions with a prescribed total liquid
contents are admissible equilibria, while, under the influence of gravity,
the only equilibrium distribution in a connected container consists in two pure phases separated by
one plane interface perpendicular to the gravity force.
\vspace{.4cm}



\vspace{1mm}

\section{Introduction}\label{intr}

In this paper we continue the discussion started in
\cite{krsbottle1}, extending the model to the case in which also gravity
effects are considered during the phase transition process.
We derive a model for solid-liquid phase transition of a medium inside
a rigid container.
The main goal is to give a qualitative and quantitative
description of the interaction between volume, pressure, phase, and temperature
changes in the situation that the specific volume of the solid phase exceeds
the specific volume of the liquid phase.
We observe, in particular, that the solidification may take place at a temperature slightly above
the critical temperature $\theta_c$. The overheating is due to the fact that the pressure decreases
from the bottom to the top. A quantitative description of this phenomenon is given in here by
means of the so-called Clausius-Clapeyron formula (cf.~\eqref{cce}).

There is an abundant classical literature on the study of phase transition processes, see e.g.
the monographs \cite{bs}, \cite{fremond}, \cite{visintin} and the references therein.
In particular, in \cite{fr1}, the authors proposed to interpret a phase transition process in terms of a
balance equation for macroscopic motions, and to include the possibility of voids, while
the microscopic approach has been pursued in \cite{fr2} in the case of two different
densities $\varrho_1$ and $\varrho_2$ for the two substances undergoing phase transitions.
Let us, however, refer to the Introduction of \cite{krsbottle1} for a more detail description of
the previous works in the literature on this topic.

The forces occurring as a result of solid-liquid phase transitions
in small containers are very strong, much stronger than gravity
forces. In a bottle of water of less than one meter height
for example, they differ by at least four orders of magnitude.
From the quantitative viewpoint, the gravity effects can
thus be neglected without considerable changes of the
process on finite time intervals. They have, however, a
substantial influence on the long time behavior of the evolution
system. Without gravity, we observe a temperature interval, in
which all phase distributions with a prescribed total liquid
contents are admissible equilibria (cf.~also \cite{krsbottle1}).
If even a weak gravity field
is assumed to be present, then the only equilibrium distribution
in a connected container consists in two pure phases separated by
one plane interface perpendicular to the gravity force.

Here we proceed as follows: in Section~\ref{mode}, we derive a model describing the evolution
of the process which is driven by an energy balance, a quasistatic
momentum balance, and a phase dynamics equation.
Still in Section~\ref{mode}, we verify the thermodynamic
consistency of the model, and we study the equilibria.

The well-posedness of the corresponding system of evolution equations is
proved  in Section \ref{exiuni}, while
and the study of the long-time behavior of solutions and convergence to equilibria
is proved in the last Section \ref{long}.

\section{The model}\label{mode}

As reference state, we consider a liquid substance contained in a bounded connected
bottle $\Omega \subset \real^3$ with boundary of class $C^{1,1}$. The state variables are
the absolute temperature $\theta>0$, the displacement $\bfu \in \real^3$,
and the phase variable $\chi \in [0,1]$. The value  $\chi = 0$
means solid, $\chi = 1$ means liquid, $\chi \in (0,1)$ is a mixture
of the two.\\[2mm]
We make the following modeling hypotheses.
\begin{itemize}
\item[{\bf (A1)}]
The displacements are small.
Therefore, we state the problem in {\em Lagrangian coordinates\/},
in which the mass conservation is equivalent to the condition
of a constant mass density $\varrho_0>0$.
\item[{\bf (A2)}]
The substance is compressible, and the speed of sound does not depend on the phase.
\item[{\bf (A3)}]
The evolution is slow, and we neglect shear viscosity and inertia effects.
\item[{\bf (A4)}]
We neglect shear stresses.
\end{itemize}
In agreement with {\bf (A1)}, we define the strain $\bfe$ as an element
of the space $\tens$ of symmetric tensors by the formula
\begin{equation}\label{eps}
\bfe = \nabla_s \bfu := \frac12 (\nabla \bfu + (\nabla \bfu)^T).
\end{equation}
Let $\bfde\in \tens$ denote the Kronecker tensor.
By {\bf (A4)}, the elasticity matrix $\bfA$ has the form
\begin{equation}\label{ela}
\bfA\bfe = \lambda (\bfe\bfdd\bfde)\,\bfde\,,
\end{equation}
where ``$\bfdd$'' is the canonical scalar product in $\tens$,
and $\lambda > 0$ is the Lam\'e constant (or {\em bulk elasticity modulus}),
which we assume to be independent of $\chi$ by virtue of {\bf (A2)}.
Note that $\lambda$ is related to the speed of sound $v_0$ by the formula
$v_0 = \sqrt{\lambda/\varrho_0}$.\\[2mm]
We want to model the situation where the specific volume $V_{solid}$
of the solid phase is larger than the specific volume $V_{liquid}$ of the
liquid phase. Considering the liquid phase as the reference state, we
{introduce} the dimensionless phase expansion coefficient
$\alpha = (V_{solid}-V_{liquid})/V_{liquid} > 0$,
and we define the phase expansion strain $\tbfe$ by
\begin{equation}\label{tbfe}
\tbfe(\chi) = \frac{\alpha}{3}(1-\chi) \bfde\,.
\end{equation}
We fix positive constants $c_0$ (specific heat), $L_0$ (latent heat),
$\theta_c$ (freezing point at standard atmospheric pressure), $\beta$ (thermal expansion coefficient), and consider the
specific free energy $f$ in the form
\begin{eqnarray}\label{free}
f &=& c_0\theta\Big(1-\log\Big(\frac{\theta}{\theta_c}\Big)\Big) + \frac{\lambda}{2\varrho_0}
((\bfe - \tbfe(\chi))\bfdd\bfde)^2 - \frac{\beta}{\varrho_0}(\theta - \theta_c)\bfe\bfdd\bfde
\\[2mm]\nonumber
&& +\, L_0\left(\chi\left(1 - \frac{\theta}{\theta_c}\right)
+  I(\chi)\right)\,,
\end{eqnarray}
where $I$ is the indicator function of the interval $[0,1]$.
The stress tensor $\bfsi$ is decomposed into the sum $\bfsi^v + \bfsi^e$
of the viscous component  $\bfsi^v$ and elastic component $\bfsi^e$.
The state functions $\bfsi^v$,$\bfsi^e$, $s$ (specific entropy),
and $\bla{e}$ (specific internal energy) are given by the formulas
\begin{eqnarray}\label{e1}
\bfsi^v &=& \nu (\bfe_t\bfdd\bfde)\bfde\\[2mm]\label{e2a}
\bfsi^e &=& \varrho_0\frac{\partial f}{\partial \bfe} =
\left(\lambda (\bfe\bfdd\bfde - \alpha(1-\chi)) - \beta (\theta - \theta_c) \right)
\bfde \,,\\[2mm]\label{e2}
s &=& -\frac{\partial f}{\partial \theta} =
c_0\log\left(\frac{\theta}{\theta_c}\right)  + \frac{L_0}{\theta_c}\chi
+ \frac{\beta}{\varrho_0} \bfe\bfdd\bfde \,,\\[2mm]\label{e3}
e &=& f + \theta\,s = c_0\theta  + \frac{\lambda}{2\varrho_0} (\bfe\bfdd\bfde - \alpha(1-\chi))^2
+ \frac{\beta}{\varrho_0}\theta_c \bfe\bfdd\bfde + L_0(\chi + I(\chi))\,,
\end{eqnarray}
\bla{where} $\nu>0$ is the volume viscosity coefficient. The scalar quantity
\begin{equation}\label{press}
p := -\nu \bfe_t\bfdd\bfde - \lambda (\bfe\bfdd\bfde - \alpha(1-\chi))
+ \beta (\theta - \theta_c)
\end{equation}
is the {\em pressure\/}, and the stress has the form $\bfsi = -p\,\bfde$.
The process is governed by the balance equations
\begin{eqnarray}
- \dive \bfsi &=& \bff_{vol} \hspace{13mm}  \mbox{(mechanical equilibrium)}\label{bal1}\\[2mm]
\varrho_0 e_t + \dive \bfq &=& \bfsi\bfdd\bfe_t \hspace{10mm}  \mbox{(energy balance)}\label{bal2} \\[2mm]
-\gamma_0\chi_t &\in& \partial_\chi f \hspace{12.5mm}  \mbox{(phase relaxation law)}\label{bal3}
\end{eqnarray}
where $\gamma_0$ is the phase relaxation coefficient, $\partial_\chi$ is the partial subdifferential
with respect to $\chi$,
$\bff_{vol}$ is a given volume force density (the gravity force)
\begin{equation}\label{grav}
\bff_{vol} = -\varrho_0 g\,\bfde_3\,,
\end{equation}
with standard gravity $g$ and vector $\bfde_3 = (0,0,1)$,
and $\bfq$ is the heat flux vector that we assume in the form
\begin{equation}\label{flux}
\bfq = -\kappa\nabla\theta
\end{equation}
with a constant heat conductivity $\kappa > 0$.
The equilibrium equation (\ref{bal1}) can be rewritten in the form $\nabla p = -\varrho_0 g$, hence
\begin{equation}\label{bala1}
p(x,t) = P(t) - \varrho_0 g\,x_3
\end{equation}
with a function $P$ of time only, which is to be determined.
On $\partial \Omega$, we assume boundary conditions in the form
\begin{eqnarray}\label{bcue}
\bfu &=& 0
\\[2mm]\label{bctheta}
\bfq\cdot \bfn &=& h(x) (\theta - \theta_\Gamma)\,,
\end{eqnarray}
with a given positive measurable function
$h$ (heat transfer coefficient), and a constant $\theta_\Gamma >0$ (external
temperature). Identity (\ref{bcue}) means that the boundary is rigid. Other
possibilities (elastic or elastoplastic boundary response) have been considered
in another context (\cite{krsbottle1,krsplast}).

By Gauss' Theorem, we have
\begin{equation}\label{bcu2}
\io \dive \bfu(x,t)\,\dd x \ = \ 0
\end{equation}
We have $\bfe\bfdd\bfde = \dive \bfu$. Using \eqref{press}, we write the mechanical equilibrium equation (\ref{bala1}) as
\begin{equation}\label{equau}
\nu \dive \bfu_t + \lambda (\dive \bfu - \alpha(1-\chi)) - \beta
(\theta - \theta_c) + P(t) - \varrho_0 g\,x_3= 0\,.
\end{equation}
Integrating over $\Omega$ and using (\ref{bcu2}) we obtain
\begin{equation}\label{equau2}
P(t) = \frac{\alpha\lambda}{|\Omega|} \io(1-\chi)\dd x' + \frac{\beta}{|\Omega|}
\io(\theta - \theta_c)\dd x' + \frac{\varrho_0 g}{|\Omega|}\io x_3'\dd x'\,.
\end{equation}
We see that in liquid ($\chi = 1$) and at temperature $\theta = \theta_c$,
the pressure $p(x,t)$ vanishes on the ``midsurface''
of $\Omega$ given by the equation $x_3 |\Omega| = \io x_3'\dd x'$.
Hence, $p(x,t)$ can be interpreted as the difference between the absolute pressure
and the standard pressure. This difference is higher below and lower above the
midsurface.

Eq.~(\ref{equau2}) enables us to eliminate $P(t)$ and rewrite (\ref{equau}) in the form
\begin{eqnarray}\nonumber
&&\hspace{-18mm}\nu \dive \bfu_t + \lambda (\dive \bfu - \alpha(1-\chi)) - \beta
(\theta - \theta_c)
+ \frac{1}{|\Omega|} \io \big(\alpha\lambda(1-\chi)+\beta(\theta - \theta_c)\big)
\dd x'\\[1mm]\label{equaup}
&=& \varrho_0 g\left(x_3 - \frac{1}{|\Omega|}\io x_3'\dd x' \right)\,.
\end{eqnarray}
As a consequence of (\ref{free}), the energy balance and the phase relaxation equation
in (\ref{bal2})--(\ref{bal3}) have the form
\begin{eqnarray}\label{equ1}
\varrho_0 c_0\theta_t - \kappa\Delta\theta &=& \nu (\dive \bfu_t)^2 - \beta\theta \dive \bfu_t -
\big(\alpha\lambda(\dive \bfu - \alpha(1-\chi)) + \varrho_0 L_0\big)\chi_t\,,\qquad\\ [2mm]\label{equ4}
 -\varrho_0\gamma_0 \chi_t &\in& \alpha\lambda
(\dive\bfu - \alpha(1-\chi)) + \varrho_0 L_0 \left(1 - \frac{\theta}{\theta_c}
+ \partial I(\chi)\right)\,,
\end{eqnarray}
where $\partial$ denotes the subdifferential. For simplicity, we now set
\begin{equation}\label{zero}
c := \varrho_0 c_0\,, \quad \gamma := \varrho_0\gamma_0\,, \quad L := \varrho_0 L_0\,.
\end{equation}
The system now completely decouples. For the unknown absolute temperature $\theta$,
local relative volume increment $U = \dive\bfu$, and liquid fraction $\chi$,
we have the evolution system (note that mathematically, $\partial I(\chi)$
is the same as $L \partial I(\chi)$)
\begin{eqnarray}\label{sys2}
c\theta_t - \kappa\Delta\theta &=& \nu U_t^2 - \beta\theta U_t -
\big(\alpha\lambda(U - \alpha(1-\chi)) + L\big)\chi_t\,,\\ [2mm]\label{sys1}
\nu U_t + \lambda U &=&
\alpha\lambda (1-\chi) + \beta (\theta - \theta_c)
+ \varrho_0 g\left(x_3  - \frac{1}{|\Omega|}\io x_3'\dd x'\right)\\ [2mm]\nonumber
&&- \frac{1}{|\Omega|} \io \big(\alpha\lambda(1-\chi)+\beta(\theta - \theta_c)\big)
\dd x'\,,\\[2mm]\label{sys3}
-\gamma \chi_t &\in& \alpha\lambda (U - \alpha(1-\chi))
+ L \left(1 - \frac{\theta}{\theta_c}\right) + \partial I(\chi)\,,
\end{eqnarray}
with boundary condition (\ref{bctheta}),
(\ref{flux}), that is,
\begin{equation}\label{bcfin}
\kappa\nabla\theta\cdot\bfn + h(x)(\theta-\theta_\Gamma) = 0\,.
\end{equation}
We then find $\bfu$ as a solution to the equation $\dive \bfu = U$ in $\Omega$, $\bfu = 0$ on $\partial\Omega$.
It is indeed not unique, and due to our hypotheses {\bf (A3)}, {\bf (A4)}, we lose any control on
possible volume preserving turbulences. This, however, has no influence
on the system (\ref{sys2})--(\ref{sys3}), which is the subject of our interest here.

Let us describe the set of all possible stationary states. It follows from (\ref{sys2}) and
(\ref{bcfin}) that the only temperature equilibrium is $\theta = \theta_\Gamma$.
The equilibrium values $\chi_\infty$ and $U_\infty$ satisfy the system
\begin{eqnarray}\label{sys1eq}
\lambda (U_\infty - \alpha (1-\chi_\infty)) &=&- \frac{\alpha\lambda}{|\Omega|}
\io (1-\chi_\infty)\dd x' + \varrho_0 g\left(x_3  - \frac{1}{|\Omega|}\io x_3'\dd x'\right),\qquad
\\[2mm]\label{sys3eq}
-\lambda (U_\infty - \alpha(1-\chi_\infty)) &\in&
\frac{L}{\alpha} \left(1 - \frac{\theta_\Gamma}{\theta_c}\right) + \partial I(\chi_\infty)\,,
\end{eqnarray}
almost everywhere in $\Omega$, that is,
\begin{equation}\label{equi1}
\frac{L}{\alpha} \left(\frac{\theta_\Gamma}{\theta_c} - 1\right)
+ \frac{\alpha\lambda}{|\Omega|} \io (1-\chi_\infty)\dd x'
 - \varrho_0 g\left(x_3  - \frac{1}{|\Omega|}\io x_3'\dd x'\right)\in\partial I(\chi_\infty(x))\,.
\end{equation}
We claim that unlike in the case without gravity, (\ref{equi1}) determines the
equilibria uniquely.
The set $\Omega$ is connected. We can therefore define
$(a,b)\subset \real$ as the maximal interval such that
$\Omega \cap (\real^2\times \{x_3\}) \ne \emptyset$ for $x_3 \in (a,b)$,
and set
\begin{equation}\label{defm}
m = \frac{1}{|\Omega|}\io x_3'\dd x' \in (a,b)\,, \quad \ell = b - a\,.
\end{equation}
We introduce two dimensionless constants
$$
d = \frac{\alpha^2\lambda}{L}\,, \qquad G_0 = \frac{\alpha\lambda}{\varrho_0 g \ell}\,.
$$
For water, we have for instance $d \approx 0.055$, $G_0 \approx 2\cdot 10^4/\ell$.
Eq.~(\ref{equi1}) then reads
\begin{equation}\label{equi2a}
\frac{1}{d} \left(\frac{\theta_\Gamma}{\theta_c} - 1\right)
+ \frac{1}{|\Omega|} \io (1-\chi_\infty)\dd x' + \frac{1}{G_0\ell}
(m  - x_3)\in\partial I(\chi_\infty(x))\,.
\end{equation}
The quantity
\begin{equation}\label{equi4}
Z := \frac{1}{d} \left(\frac{\theta_\Gamma}{\theta_c} - 1\right)
+ \frac{1}{|\Omega|} \io (1-\chi_\infty)\dd x'
\end{equation}
is independent of $x$, so that the left hand side of (\ref{equi2a}) is positive for
$x_3 < m + G_0 \ell Z$ and negative for $x_3 > m + G_0 \ell Z$. By definition of the
subdifferential, we necessarily have
\begin{equation}\label{equi3a}
\chi_\infty(x) = \left\{
\begin{array}{ll}
1 & \pif\ x_3 < m + G_0\ell Z\,,\\[2mm]
0 & \pif\ x_3 > m + G_0\ell Z\,.
\end{array}
\right.
\end{equation}
Let $\Omega(r)$ denote the set $\{x \in \Omega: x_3 > r\}$ for $r \in \real$.
Eq.~(\ref{equi3a}) states that
the set $\Omega(m +  G_0\ell Z)$ corresponds to the solid domain.
We have $|\Omega(r)| = 0$ for $r\ge b$, $|\Omega(r)| = |\Omega|$ for $r\le a$, and
\begin{equation}\label{equi5a}
\frac{1}{|\Omega|} \io (1-\chi_\infty)\dd x' = F(Z):=
\frac{|\Omega(m +  G_0\ell Z)|}{|\Omega|}\,.
\end{equation}
We easily identify $Z$ as the only solution to the equation
\begin{equation}\label{equi6}
Z = \frac{1}{d} \left(\frac{\theta_\Gamma}{\theta_c} - 1\right)
+ F(Z)\,,
\end{equation}
since $F$ is nonincreasing. We see that one of the following three cases necessarily occurs:
\begin{itemize}
\item[(i)] $Z \le (a-m)/(G_0\ell)$, $F(Z) = 1$. Then $\chi_\infty(x) = 0$ a.e. in $\Omega$
and we have pure solid with temperatures
$$\theta_\Gamma = \theta_c (1 + d(Z-1))
\le \theta_c \left(1 - {d} \left(1 + \frac{m-a}{G_0\ell}\right)\right);
$$
\item[(ii)] $Z \ge (b-m)/(G_0\ell)$, $F(Z) = 0$. Then $\chi_\infty(x) = 1$ a.e. in $\Omega$
and we have pure liquid with temperatures
$$\theta_\Gamma = \theta_c (1 + dZ)
\ge \theta_c \left(1 + d\, \frac{b-m}{G_0\ell}\right);
$$
\item[(iii)]  $(a-m)/(G_0\ell) < Z < (b-m)/(G_0\ell)$, $0<F(Z) < 1$. Then
$\chi_\infty(x) = 0$ a.e. in $\Omega(m +  G_0\ell Z)$,
$\chi_\infty(x) = 1$ a.e. in $\Omega\setminus \Omega(m +  G_0\ell Z)$, and
$$
\theta_c \left(1 - {d} \left(1 + \frac{m-a}{G_0\ell}\right)\right) <
\theta_\Gamma <  \theta_c \left(1 + d\, \frac{b-m}{G_0\ell}\right).
$$
\end{itemize}
We observe that solidification may take place at temperatures slightly above
$\theta_c$.
For water in a container of $\ell =50\,cm$ height, the relative size of the
``overheated ice temperature domain''
is smaller than $d/G_0 \approx 1.4\cdot 10^{-6}$, hence it is far beyond
the standard measurement accuracy. The overheating is due to the fact
the pressure decreases from the bottom to the top,
as pointed out after formula (\ref{equau2}). A quantitative characterization of this
phenomenon is given by the so-called Clausius-Clapeyron equation,
which relates the freezing temperature with the pressure. It can be derived here as follows.
The equilibrium relative pressure $p_\infty$ depends only on $x_3$, and is given,
by virtue of (\ref{bala1}) and (\ref{equau2}), by the formula
\begin{equation}\label{equau3}
p_\infty(x_3) = \frac{\alpha\lambda}{|\Omega|} \io(1-\chi_\infty)\dd x' + \beta
(\theta_\Gamma - \theta_c) + \frac{\varrho_0 g}{|\Omega|}\io x_3'\dd x'
- \varrho_0 g x_3\,.
\end{equation}
Using (\ref{equi1}), we obtain
\begin{equation}\label{equi1cc}
\frac{L}{\alpha} \left(\frac{\theta_\Gamma}{\theta_c} - 1\right)
-\beta (\theta_\Gamma - \theta_c) + p_\infty(x_3)
\in\partial I(\chi_\infty(x))\,.
\end{equation}
The phase interface at temperature $\theta_\Gamma$ is located at level $x_3$ if the right hand
side of (\ref{equi1cc}) vanishes.
Setting $L_\beta = L_0 - \alpha\beta\theta_c/\varrho_0$, we thus obtain the Clausius-Clapeyron
condition for phase transition in the form of \cite[Equation (288)]{ye}, namely
\begin{equation}\label{cce}
\frac{p_\infty(x_3)}{\theta_\Gamma - \theta_c} = - \frac{\varrho_0 L_\beta}{\alpha\theta_c}
= \frac{L_\beta}{(V_{liquid} - V_{solid})\theta_c}\,.
\end{equation}

In terms of the new variables $\theta, U, \chi$, the energy $e$ and entropy $s$ can be written~as
\begin{eqnarray}\label{energy}
e &=& c_0\theta + \frac{\lambda}{2\varrho_0}(U - \alpha(1-\chi))^2
+ \frac{\beta}{\varrho_0} \theta_c U
+ L_0(\chi + I(\chi))\,,\\[2mm]\label{entropy}
s &=& c_0\log\left(\frac{\theta}{\theta_c}\right)  + \frac{L_0}{\theta_c}\chi
+ \frac{\beta}{\varrho_0} U \,,
\end{eqnarray}
\bla{and} the energy and entropy balance equations now read
\begin{eqnarray}\label{princ1}
\frac{\dd}{\dd t} \io (\varrho_0 e(x,t) - \varrho_0 g x_3 U(x,t))\dd x
&=& \ipo h(x)(\theta_\Gamma - \theta)\dd s(x)\,,\\[2mm]\label{princ2}
\varrho_0 s_t + \dive \frac{\bfq}{\theta} &=& \frac{\kappa |\nabla\theta|^2}{\theta^2}
+ \frac{\gamma}{\theta} \chi_t^2 + \frac{\nu}{\theta}U_t^2 \ \ge\ 0\,,\\[2mm]\label{princ3}
\frac{\dd}{\dd t}\io \varrho_0 s(x,t)\dd x &=& \ipo \frac{h(x)}{\theta}(\theta_\Gamma - \theta)
\dd s(x) \\[1mm]\nonumber
&& +\, \io\left(\frac{\kappa |\nabla\theta|^2}{\theta^2}
+ \frac{\gamma}{\theta} \chi_t^2 + \frac{\nu}{\theta}U_t^2  \right)\dd x\,.
\end{eqnarray}
The entropy balance (\ref{princ2}) says that the entropy production
on the right hand side is nonnegative in agreement with the second principle
of thermodynamics. The system is not closed, and the energy supply through
the boundary is given by the right hand side of (\ref{princ1}).

We prescribe the initial conditions
\begin{eqnarray}\label{ini2}
\theta(x,0) &=& \theta^0(x)\\ \label{ini1}
U(x,0) &=& U^0(x)\\ \label{ini3}
\chi(x,0) &=& \chi^0(x)
\end{eqnarray}
for $x \in \Omega$, such that $\io U_0(x)\dd x = 0$, and compute from
(\ref{energy})--(\ref{entropy}) the corresponding
initial values $e^0$ and $s^0$ for specific energy and entropy, respectively.
Let $E^0 = \io \varrho_0 e^0\dd x$, $S^0 = \io \varrho_0 s^0 \dd x$
denote the total initial energy and entropy, respectively.
{}From the energy end entropy balance equations (\ref{princ1}), (\ref{princ3})
and using the condition $\io U\dd x = 0$, we derive
the following crucial (formal for the moment) balance equation for the ``extended'' energy
$\varrho_0 (e - \theta_\Gamma s)$:
\begin{eqnarray}\label{crucial}
&&\hspace{-16mm}\io\left(c\theta + \frac{\lambda}{2}(U - \alpha(1-\chi))^2
+ L\chi - \varrho_0 g x_3 U \right)(x,t)\dd x\\ \nonumber
&&+\, \theta_\Gamma \int_0^t\io \left(\frac{\kappa |\nabla\theta|^2}{\theta^2}
+ \frac{\gamma}{\theta} \chi_t^2 + \frac{\nu}{\theta}U_t^2  \right)(x,\tau)\dd x \dd\tau\\ \nonumber
&&+\, \int_0^t \ipo\frac{h(x)}{\theta}(\theta_\Gamma - \theta)^2(x,\tau)\dd s(x) \dd\tau\\ \nonumber
&=& E^0 - \theta_\Gamma S^0 -\varrho_0 g \io x_3 U^0(x) \dd x
+ \theta_\Gamma\io\left(c\log\left(\frac{\theta}{\theta_c}\right)  +
\frac{L}{\theta_c}\chi \right)(x,t)\dd x\,.
\end{eqnarray}
We have $\log(\theta/\theta_c) = \log(\theta/2\theta_\Gamma) - \log(\theta_c/2\theta_\Gamma)
\le (\theta/2\theta_\Gamma) - 1 - \log(\theta_c/2\theta_\Gamma)$, hence there exists
a constant $C>0$ independent of $t$ such that for all $t>0$ we have
\begin{eqnarray}\label{esti1}
&&\hspace{-16mm}\io\left(\theta + U^2\right)(x,t)\dd x
+ \int_0^t\io \left(\frac{|\nabla\theta|^2}{\theta^2}
+ \frac{\chi_t^2}{\theta}  + \frac{U_t^2}{\theta}  \right)(x,\tau)\dd x \dd\tau\\ \nonumber
&&+\, \int_0^t \ipo\frac{h(x)}{\theta}(\theta_\Gamma - \theta)^2(x,\tau)\dd s(x) \dd\tau \ \le \ C\,.
\end{eqnarray}

\section{Existence and uniqueness of solutions}\label{exiuni}

We construct the solution of (\ref{sys1})--(\ref{sys3}) by the Banach
contraction argument. The method of proof is independent
of the actual values of the material constants, and we choose for simplicity
\begin{equation}\label{const}
L= 2,\ \ g=c= \varrho_0 =\theta_c =\alpha = \beta= \gamma = \kappa = \lambda = \nu = 1\,.
\end{equation}
System (\ref{sys2})--(\ref{sys3}) with boundary condition (\ref{bctheta}) then reads
\begin{eqnarray}\label{nequ1}
\io\theta_t w(x)\dd x + \io \nabla\theta \cdot \nabla w(x) \dd x &=&
\io\left(U_t^2 - \theta U_t -\Big(U +\chi +1\Big)\chi_t\right) w(x)\dd x\qquad \\ \nonumber
&&-\,\ipo h(x)(\theta - \theta_\Gamma) w(x) \dd s(x) \,,\\ [2mm]\label{nequ2}
U_t + U &=&-\chi + \theta
+ \left(x_3  - \frac{1}{|\Omega|}\io x_3'\dd x'\right)\\ [2mm]\nonumber
&&- \frac{1}{|\Omega|} \io \big(-\chi+\theta\big)
\dd x'\,,\\ [2mm]\label{nequ3}
\chi_t + U +\chi + \partial I(\chi) &\ni& 2\theta - 1\,,
\end{eqnarray}
where (\ref{nequ1}) is to be satisfied for all test functions $w \in W^{1,2}(\Omega)$
and a.e. $t>0$, while (\ref{nequ2})--(\ref{nequ3}) are supposed to hold a.e. in
$\Omega_\infty := \Omega\times (0,\infty)$.

In this section we prove the following existence and uniqueness result.
\begin{theorem}\label{main}
Let $0< \theta_* \le \theta_\Gamma \le \theta^*$ and $p_0 \in \real$
be given constants, and let the data satisfy the conditions
$$
\barr{rcllcccll}
\theta^0 &\in& W^{1,2}(\Omega)\cap L^\infty(\Omega)\,,
\quad &\theta_* &\le& \theta^0(x)&\le& \theta^*\quad &\mbox{a.e.},\\
U^0, \chi^0 &\in& L^\infty(\Omega)\,,\quad \io U^0(x)\dd x = 0\,, &0 &\le& \chi^0(x) &\le& 1 \quad &\mbox{a.e.}
\earr
$$
Then there exists a unique solution $(\theta,U,\chi)$ to (\ref{nequ1})--(\ref{nequ3}),
(\ref{ini2})--(\ref{ini3}), such that
$\theta>0$ a.e., $\chi \in [0,1]$ a.e.,
$U, U_t,\chi_t, \theta, 1/\theta \in L^\infty(\Omega_\infty)$,
$\theta_t, \Delta\theta \in L^2(\Omega_\infty)$,
and $\nabla\theta \in L^\infty(0,T;L^2(\Omega)) \cap L^2(\Omega_\infty)$.
\end{theorem}

\begin{remark}\label{bc}
\rm For existence and uniqueness alone, we might allow the external temperature
$\theta_\Gamma$ to depend on $x$ and $t$, and assume only that it belongs to the space
$W^{1,2}_{{\rm loc}} (0,\infty;L^2(\partial\Omega))
\cap L^\infty_{{\rm loc}}(\partial\Omega\times(0,\infty))$.
For the global bounds, the assumption that $\theta_\Gamma$ be constant plays a
substantial role.
\end{remark}

The proof of Theorem \ref{main} will be carried out in the following subsections.
Notice first that the term $U_t^2 - \theta U_t - (U +\chi +1)\chi_t$ on the right hand
side of (\ref{nequ1}) can be rewritten alternatively, using \eqref{nequ3} and \eqref{nequ2}, as
\begin{eqnarray}\label{rhs}
U_t^2 - \theta U_t - (U +\chi +1)\chi_t &=& U_t^2 - \theta U_t + \chi_t^2 - 2 \theta\chi_t
\\[2mm]\nonumber
&=& -(\chi+U)U_t
+\left(x_3  - \frac{1}{|\Omega|}\io x_3'\dd x'\right)U_t\\[2mm]\nonumber
&&-\left( \frac{1}{|\Omega|} \io \big(-\chi+\theta\big)\dd x' \right)U_t - (U +\chi +1)\chi_t\,,
\end{eqnarray}
We now fix some constant $R>0$ and construct the solution for the truncated system
\begin{eqnarray}\label{trequ1}
\io\theta_t w(x)\dd x + \io \nabla\theta \cdot \nabla w(x) \dd x &=&
\io\left(U_t^2 +\chi_t^2 -  Q_R(\theta) (U_t+ 2 \chi_t)\right) w(x)\dd x\qquad \\ \nonumber
&&\hspace{-15mm}-\,\ipo h(x)(\theta - \theta_\Gamma) w(x)\dd s(x)\quad \forall w \in W^{1,2}(\Omega)\,,
\\ [2mm]\label{trequ2}
U_t + U &=&-\chi + Q_R(\theta)
+ \left(x_3  - \frac{1}{|\Omega|}\io x_3'\dd x'\right)\\ [2mm]\nonumber
&&- \frac{1}{|\Omega|} \io \big(-\chi+Q_R(\theta)\big)
\dd x'\,,\\ [2mm]\label{trequ3}
 \chi_t + U +\chi + \partial I(\chi) &\ni& 2 Q_R(\theta) - 1
\end{eqnarray}
first in a bounded domain $\Omega_T := \Omega\times (0,T)$ for any given $T>0$,
where $Q_R$ is the cutoff function $Q_R(z) = \min\{z^+, R\}$. We then derive upper and lower
bounds for $\theta$ independent of $R$ and $T$, so that the local solution of
(\ref{trequ1})--(\ref{trequ3}) is also a global solution of (\ref{nequ1})--(\ref{nequ3})
if $R$ is sufficiently large.

\subsection{A gradient flow}\label{grad}

System (\ref{sys1})--(\ref{sys3}) can be considered as a gradient flow
similarly as in the case without gravity and with an elastic boundary,
see \cite[Section~4.1]{krsbottle1}. Set
\begin{eqnarray}\label{aseg4}
v &=& \vect{ U}{\chi}\,,\\[2mm]\label{aseg5}
 \quad \psi(v) &=& \io \left( \frac{1}{2} (U -(1-\chi))^2
+ 2\chi \left(1-\theta_\Gamma\right) - x_3 U
+ I(\chi)\right)\dd x\qquad \\ \nonumber
&&+\, \frac{1}{|\Omega|} \io U \dd x \io (1-\chi+x_3)\dd x
 + C_\psi\,,\\[2mm]\label{aseg6}
f &=& \vect{(\theta - \theta_\Gamma) - 1/|\Omega| \io (\theta - \theta_\Gamma)
\dd x}{2(\theta - \theta_\Gamma)}
\end{eqnarray}
with a constant $C_\psi \ge 0$ to ensure that $\psi(v) \ge 0$. Every solution $U$
of (\ref{sys1}) necessarily satisfies the condition $\io U \dd x = 0$, hence
(\ref{sys3}) can be written in the form
\begin{equation}\label{aseg7}
\chi_t + U+\chi+1 -2 \theta- \frac{1}{|\Omega|}
\io U \dd x + \partial I(\chi) \ni 0\,,
\end{equation}
and (\ref{sys1}), (\ref{aseg7}) are in turn equivalent to the gradient flow
\begin{equation}\label{aseg8}
\dot v + \partial \psi(v) \ni f
\end{equation}
in $L^2(\Omega) \times L^2(\Omega)$. We have used the obvious identity
$\theta - 1 - 1/|\Omega| \io (\theta - 1)
\dd x = (\theta - \theta_\Gamma) - 1/|\Omega| \io (\theta - \theta_\Gamma)
\dd x$.
We state here the following Lemma, whose proof can be found in \cite[Lemma~4.3]{krsbottle1}.

\begin{lemma}\label{lg1}
Let $f, \dot f$ belong to $L^2(0,\infty;H)$. Then $\lim_{t \to \infty} \dot v(t) = 0$.
\end{lemma}

We apply the above result to the case $H=L^2(\Omega)\times L^2(\Omega)$, and
$v$, $f$, $\psi$ as above
and we see that Eqs.~(\ref{trequ2})--(\ref{trequ3}) with $\theta$ replaced by $\hat\theta$
can be equivalently written as a gradient flow (\ref{aseg4})--(\ref{aseg8}).
For its solutions, we prove the following result.

\begin{proposition}\label{pg1}
Let the hypotheses of Theorem \ref{main} hold, and let a function
$\hat\theta \in L^2_{{\rm loc}}(0,\infty; L^2(\Omega))$ be given.
Let $(U,\chi)$ be the solution of
(\ref{aseg4})--(\ref{aseg8}). Then there exists a constant $C_0$, independent
of $x, t$ and $R$, such that a.e. in $\Omega_\infty$ we have
\begin{equation}\label{eg8bis}
|U(x,t)| + |U_t(x,t)| + |\chi_t(x,t)| \le C_0(1+R)\,.
\end{equation}
Let furthermore $\hat\theta_1, \hat\theta_2 \in L^2_{{\rm loc}}(0,\infty; L^2(\Omega))$
be two functions, and let $(U_1,\chi_1), (U_2,\chi_2)$ be the corresponding solutions of
(\ref{aseg4})--(\ref{aseg8}). Then the differences
$\hat\theta_d = \hat\theta_1 - \hat\theta_2$, $U_d = U_1-U_2$, $\chi_d =\chi_1-\chi_2$
satisfy for every $t \ge 0$ and a.e. $x \in \Omega$ the inequality
\begin{equation}\label{eg9}
\int_0^t (|(U_d)_t| + |(\chi_d)_t|)(x,\tau) \dd \tau \le C_0(1+t)
\int_0^t\left(|\hat\theta_d(x,\tau)| + t |\hat\theta_d(\tau)|_2\right)\dd \tau\,,
\end{equation}
where the symbol $|\cdot|_2$ stands for the norm in $L^2(\Omega)$.
\end{proposition}

In what follows,
we denote by $C_1, C_2, \dots$ any constant independent of $x, t$ and~$R$.

\bpf{Proof}
We rewrite (\ref{aseg4})--(\ref{aseg8}) as two scalar gradient
flows
\begin{eqnarray}\label{eg10}
U_t + \partial \psi_1(U) &=& a\,,\\[1mm]\label{eg11}
\chi_t + \partial \psi_2(\chi) &\ni& b\,,
\end{eqnarray}
where $\psi_1(U) = \frac12 U^2$, $\psi_2 = \frac12 \chi^2 + I(\chi)$,
$a = Q_R(\hat\theta) -\chi - \frac{1}{|\Omega|}\int_\Omega(Q_R(\hat\theta) -\chi)\dd x
+ \left(x_3  - \frac{1}{|\Omega|}\io x_3'\dd x'\right)$, $b = 2 Q_R(\hat\theta)-1-U$.
The bounds (\ref{eg8bis}) are obvious. To prove (\ref{eg9}), we
consider two different inputs. As above, we denote the differences $\{\}_1-\{\}_2$
by $\{\}_d$ for all symbols $\{\}$.
By \cite[Theorem 1.12]{hys}, we have for all $t>0$ and a.e. $x \in \Omega$ that
 \begin{equation}\label{eg12}
\int_0^t (|(U_d)_t| + |(\chi_d)_t|)(x,\tau) \dd \tau \le 2\int_0^t
(|a_d| + |b_d|)(x,\tau) \dd \tau\,.
 \end{equation}
We multiply the difference of (\ref{eg10}) by $U_d$, the difference of (\ref{eg11})
by $\chi_d$, and sum them up to obtain that
\begin{equation}\label{eg13}
(U_d)_t U_d + (\chi_d)_t\chi_d + (U_d +\chi_d)^2
\le |\hat\theta_d|(|U_d| + 2 |\chi_d|)-\left(\frac{1}{|\Omega|}
\int_\Omega(Q_R(\hat\theta)_d-\chi_d)\dd x\right)U_d \quad\mae
\end{equation}
We first integrate (\ref{eg13}) over $\Omega$. Using the symbol $|\cdot|_2$ for the norm
in $L^2(\Omega)$, we get for a.e. $t >0$ that
\begin{equation}\label{eg14}
\frac12 \frac{\dd}{\dd t} \left(|U_d|_2^2 + |\chi_d|_2^2\right)  \le
|\hat\theta_d|_2(|U_d|_2 +2 |\chi_d|_2) \le
\sqrt{5}|\hat\theta_d|_2 \left(|U_d|_2^2 +|\chi_d|_2^2\right)^{1/2}\,.
\end{equation}
Hence, $\frac{\dd}{\dd t} (|U_d|_2^2 + |\chi_d|_2^2)^{1/2} \le \sqrt{5}|\hat\theta_d|_2$
a.e., and integrating over $t$, we find that
 \begin{equation}\label{eg15}
\left(|U_d|_2^2 +|\chi_d|_2^2\right)^{1/2}(t)
\le \sqrt{5}\int_0^t |\hat\theta_d(\tau)|_2\dd\tau\,.
 \end{equation}
Using again (\ref{eg13}), we find for a.e. $(x,t)\in \Omega_\infty$ the inequality
\begin{align}\nonumber
\frac12 \frac{\partial}{\partial t} \left(U_d^2 + \chi_d^2\right)(x,t)  \le \ &|\hat\theta_d|(x,t)(|U_d| +2 |\chi_d|)(x,t)\\
\label{eg17}
&+C_1|U_d|(x,t)\left(|\hat\theta_d(t)|_2+\int_0^t|\hat\theta_d(x,\tau)|_2\dd\tau\right)\,.
\end{align}
Hence, we get
\begin{equation}\label{eg18}
\frac12 \frac{\partial}{\partial t} \left(U_d^2 + \chi_d^2\right)(x,t)
\le C_2 (U_d^2+\chi_d^2)^{1/2}(x,t)\left(|\hat\theta_d|(x,t)
+|\hat\theta_d(t)|_2 +\int_0^t|\hat\theta_d(\tau)|_2\dd\tau\right)\,,
\end{equation}
which in turn implies that
\begin{equation}\label{eg19}
\left(|U_d|^2 + |\chi_d|^2\right)^{1/2} (x,t)  \le
C_3\left(\int_0^t|\hat\theta_d(x,\tau)|\dd\tau + (1+t) \int_0^t|\hat\theta_d(\tau)|_2\dd \tau\right)
\quad \mbox{a.e.}
\end{equation}
This enables us to estimate the right hand side of (\ref{eg12}) and obtain the bound
 \begin{equation}\label{eg20}
\int_0^t (|(U_d)_t| + |(\chi_d)_t|)(x,\tau)\dd \tau \le C_4\int_0^t
\left((1+t)|\hat\theta_d|(x,\tau) + t(1+t) |\hat\theta_d(\tau)|_2\right)\dd \tau
 \end{equation}
for a.e. $x\in \Omega$ and all $t \ge 0$. This completes the proof.
\epf

\subsection{Existence of solutions for the truncated problem}\label{exitrunc}

We construct the solution of (\ref{trequ1})--(\ref{trequ3}) for every $R>0$ by the Banach
contraction argument on a fixed time interval $(0,T)$.

\begin{lemma}\label{le1}
Let the hypotheses of Theorem \ref{main} hold, and let $T>0$ and $R>0$ be given.
Then there exists
a unique solution $(\theta,U,\chi)$ to (\ref{trequ1})--(\ref{trequ3}),
(\ref{ini2})--(\ref{ini3}), such that $U\in W^{1,\infty}(\Omega_T)$, $\theta>0$ a.e.,
$\chi_t, \theta, 1/\theta \in L^\infty(\Omega_T)$, $\theta_t, \Delta\theta
\in L^2(\Omega_T)$, and $\nabla\theta \in L^\infty(0,T;L^2(\Omega))$.
\end{lemma}

\bpf{Proof}\ Let $\hat\theta \in L^2(\Omega_T)$ be a given function,
and consider the system
\begin{eqnarray}\label{appqu1}
\io\theta_t w(x)\dd x + \io \nabla\theta \cdot \nabla w(x) \dd x \!&\!=\!&\!\!
\io\left(U_t^2 +\chi_t^2 -  Q_R(\hat\theta) (U_t+ 2 \chi_t)\right) w(x)\dd x\qquad \quad \\ \nonumber
&&\hspace{-15mm}-\,\ipo h(x)(\theta - \theta_\Gamma) w(x)\dd s(x)\quad \forall w \in W^{1,2}(\Omega)\,,
\\ [2mm]\label{appqu2}
U_t + U &=&-\chi + Q_R(\hat\theta)
+ \left(x_3  - \frac{1}{|\Omega|}\io x_3'\dd x'\right)\\ [2mm]\nonumber
&&- \frac{1}{|\Omega|} \io \big(-\chi+Q_R(\hat\theta)\big)
\dd x'\,,\\ [2mm]\label{appqu3}
 \chi_t + U +\chi + \partial I(\chi) &\ni& 2 Q_R(\hat\theta) - 1
\end{eqnarray}
Equations (\ref{appqu2})--(\ref{appqu3}) are solved as a gradient flow problem from
Subsection \ref{grad}, while (\ref{appqu1}) is a simple linear parabolic equation for $\theta$.
Testing (\ref{appqu1}) by $\theta_t$, we obtain by Proposition \ref{pg1} that
\begin{eqnarray}\label{eg19a}
&&\int_0^T\io \theta_t^2\dd x \dd t + \supess_{t\in (0,T)}\left(
\io |\nabla\theta|^2 \dd x + \ipo h(x)(\theta - \theta_\Gamma)^2 \dd s(x)\right)
\\[2mm]\nonumber
&& \qquad \le T |\Omega| \left(C_0(1+R)(2 C_0(1+R) + 3R)\right)^2 =: M_R\,.
\end{eqnarray}
Hence, we can define the mapping that with $\hat\theta$ associates the solution $\theta$
of (\ref{appqu1})--(\ref{appqu3}) with initial conditions (\ref{ini2})--(\ref{ini3}).
We now show that it is a contraction on the set
\begin{equation}\label{eg8}
\Xi_{T,R} := \{\hat\theta \in L^2(\Omega_T): \mbox{conditions (\ref{eg8a})--(\ref{eg8d})
hold}\}\,,
\end{equation}
where
\begin{eqnarray}\label{eg8a}
&&\hspace{-14mm}\hat\theta_t \in L^2(\Omega_T)\,;\\[2mm]\label{eg8b}
&&\hspace{-14mm}\nabla\hat\theta \in L^\infty(0,T;L^2(\Omega))\,;\\[2mm]\label{eg8c}
&&\hspace{-14mm}\int_0^T\io \hat\theta_t^2 \dd x \dd t + \supess_{t\in (0,T)}\left(
\io |\nabla\hat\theta|^2 \dd x + \ipo h(x)(\hat\theta - \theta_\Gamma)^2 \dd s(x)\right)
\le M_R\,;\\[2mm]\label{eg8d}
&&\hspace{-14mm}\hat\theta(x,0) = \theta^0(x)\ \mae
\end{eqnarray}
Let $\hat\theta_1, \hat\theta_2$
be two functions in $\Xi_{T,R}$, and let $(\theta_1, U_1, \chi_1)$, $(\theta_2, U_2, \chi_2)$,
be the corresponding solutions to (\ref{appqu1})--(\ref{appqu3}) with the same initial
conditions $\theta^0, U^0,\chi^0$. We see from (\ref{eg19a}) that
$\theta_1, \theta_2$ belong to $\Xi_{T,R}$.
Integrating Eq.~(\ref{appqu1}) for $\theta_1$ and $\theta_2$ with respect to time and testing
their difference by $w = \theta_d := \theta_1-\theta_2$, we obtain, using Proposition~\ref{pg1}, that
\begin{eqnarray}\nonumber
&& \hspace{-15mm}\io\theta_d^2(x,t)\dd x +
\frac{\dd}{\dd t}\left(\io\left|\nabla\int_0^t \theta_d(x,\tau)\dd\tau\right|^2\dd x
+ \ipo h(x) \left|\int_0^t \theta_d(x,\tau)\dd\tau\right|^2\dd s(x)\right)\\[2mm]  \label{lipthe2}
 &\le& C_5(1+R) \io \left(\int_0^t (|(U_d)_t| + |(\chi_d)_t| + |\hat\theta_d|)
(x,\tau) \dd \tau\right) \theta_d(x,t) \dd x\quad \mae
\end{eqnarray}
{}From (\ref{eg9}) and Minkowski's inequality, it follows that
\begin{eqnarray*}
\left|\int_0^t (|(U_d)_t| + |(\chi_d)_t|)(\cdot,\tau) \dd \tau\right|_2
&\le& C_6 (1+t)^2\int_0^t |\hat\theta_d(\tau)|_2 \dd \tau\\[2mm]
&\le& C_6 (1+t)^2\left(t\int_0^t |\hat\theta_d(\tau)|_2^2 \dd \tau\right)^{1/2}.
\end{eqnarray*}
By Young's inequality, we rewrite (\ref{lipthe2}) as
\begin{eqnarray}\nonumber
&& \hspace{-15mm}\io\theta_d^2(x,t)\dd x +
\frac{\dd}{\dd t}\left(\io\left|\nabla\int_0^t \theta_d(x,\tau)\dd\tau\right|^2\dd x
+ \ipo h(x) \left|\int_0^t \theta_d(x,\tau)\dd\tau\right|^2\dd s(x)\right)\\[2mm]  \label{lipthe4}
 &\le& C_7(1+R^2)(1+t)^5\int_0^t |\hat\theta_d(\tau)|_2^2 \dd \tau \quad \mae
\end{eqnarray}
Set $\Theta^{2} (t) = \int_0^t|\theta_d(\tau)|_2^2 \dd\tau$,
$\hat\Theta^{2} (t) = \int_0^t|\hat\theta_d(\tau)|_2^2 \dd\tau$.
Integrating (\ref{lipthe4}) with respect to time, we obtain
\begin{equation}\label{lipthe3}
\Theta^{2}(t) \le  C_7(1+R^2) \int_0^t (1+\tau)^5 \hat\Theta^{2}(\tau)\dd\tau\,.
\end{equation}
We set $C_R:= (C_7(1+R^2)/{6})$ and introduce in $L^\infty(0,T)$ the norm
$$
\|w\|_C := \sup_{\tau \in [0,T]} \expe^{- C_R(1+\tau)^{6}}|w(\tau)|\,.
$$
Then $\|\Theta\|_C^{2} \le \frac12\|\hat\Theta\|_C^{2}$, and
hence the mapping $\hat\theta\mapsto \theta$
is a contraction in $L^2(\Omega_T)$ with respect to the norm induced by $\|\cdot\|_C$.
The set $\Xi_{T,R}$ is a closed subset of $L^2(\Omega_T)$. This implies the existence of a
fixed point $\theta \in \Xi_{T,R}$, which is indeed a solution to (\ref{trequ1})--(\ref{trequ3}).
The positive upper and lower bounds for $\theta$ follow from the maximum principle
(for the proof cf.~\cite[Section~4.2]{krsbottle1}).
This completes the proof of Lemma \ref{le1}.
\epf

\subsection{Proof of Theorem \ref{main}}\label{concluwellp}

The unique solution $(\theta,U,\chi)$ to (\ref{trequ1})--(\ref{trequ3}),
(\ref{ini2})--(\ref{ini3}) exists globally in the whole domain $\Omega_{\infty}$.
We now derive uniform bounds independent of $t$ and $R$.
Take first for instance any $R> 2\theta^*$. We know that
the solution component $\theta$ of (\ref{trequ1})--(\ref{trequ3}) remains
smaller than $R$ in a nondegenerate interval $(0,T)$ with $T > \theta^*/(C_6(1+R)^2)$.
Let $(0,T_0)$ be the maximal interval in which $\theta$ is bounded by $R$.
Then, in $(0,T_0)$, the solution given by Lemma~\ref{le1} is also a solution of the original problem
(\ref{nequ1})--(\ref{nequ3}).
Moreover, due to estimate \eqref{esti1}, we know that $\theta$ admits a bound
in $L^\infty(0,T_0;L^1(\Omega))$ independent of $R$. In order to prove that
$T_0 = +\infty$ if $R$ is sufficiently large, we need the following
variant of the Moser iteration lemma, whose proof can be found
in \cite[Prop.~4.6]{krsbottle1}.

\begin{proposition}\label{moser}
Let $\Omega \subset \real^N$ be a bounded domain with Lipschitzian boundary.
Given nonnegative functions $h \in L^1(\partial\Omega)$
and $r \in L^\infty(0,\infty; L^{q} (\Omega))$ with a fixed $q > N/2$,
$|r|_{L^\infty(0,\infty; L^q(\Omega))} =: r^*$,
an initial condition $v^0\in L^\infty(\Omega)$, and a boundary datum
$v_\Gamma \in L^\infty(\partial\Omega\times (0,\infty))$, consider the problem
\begin{eqnarray}\label{moser1}
v_t -\Delta v + v &=& r(x,t)\,{\mathcal H}[v]\qquad\hbox{a.e.~in }
\Omega \times (0,\infty)\,,\\
\label{moser1a}
\nabla v\cdot \bfn &=& - h(x)\,\left(f(x,t, v(x,t)) - v_\Gamma(x,t)\right)
\quad\hbox{a.e.~on }\partial\Omega \times (0,\infty)\,,\qquad\\
\label{moser2}
v(x,0)&=&v^0\qquad \hbox{a.e.~in }\Omega\,,
\end{eqnarray}
under the assumption that there exist positive constants
$m, H_0, C_f, V, V_\Gamma, E_0$
such that the following holds:
\begin{itemize}
\item[{\rm (i)}] The mapping ${\mathcal H}\,:\,L^\infty_{{\rm loc}}
(\Omega \times (0,\infty))\to
L^\infty_{\rm loc}(\Omega \times (0,\infty))$ satisfies for every
$v\in L^\infty_{\rm loc}(\Omega \times (0,\infty))$ and a.e.
$(x,t) \in \Omega \times (0,\infty)$ the inequality
\[
v(x,t)\,{\mathcal H}[v](x,t) \leq H_0|v(x,t)| \left(1+ |v(x,t)| + \int_0^t
\xi(t-\tau)|v(x,\tau)| \dd\tau\right),
\]
where $\xi \in W^{1,1}(0,\infty)$ is a given nonnegative function such that
\begin{equation}\label{me1a}
\dot \xi(t) \le -\xi(0)\,\xi(t) \quad \mbox{a.e.}
\end{equation}
\item[{\rm (ii)}] $f$ is a Carath\'eodory function on $\Omega \times (0,\infty) \times \real$
such that $f(x,t,v)\,v\geq C_f\,v^2$ a.e. for all $v\in \real$.
\item[{\rm (iii)}] $|v^0(x)|\leq V$ a.e.~in $\Omega$.
\item[{\rm (iv)}] $|v_\Gamma(x,t)|\leq V_\Gamma$ a.e.~on
$\partial\Omega \times (0,\infty)$.
\item[{\rm (v)}] System (\ref{moser1})--(\ref{moser2}) admits a solution\\
$v \in W^{1,2}_{{\rm loc}}(0,\infty;
(W^{1,2})'(\Omega))\cap L^2_{{\rm loc}}(0,\infty;W^{1,2}(\Omega))
\cap L^\infty_{\rm loc}(\Omega \times (0,\infty))$\\
satisfying the estimate
\[
\io|v(x,t)|\,dx\leq E_0 \quad\hbox{a.e.~in }(0,\infty)\,.
\]
\end{itemize}
Then there exists a positive constant $C^*$ depending only on
$|h|_{L^1(\partial\Omega)}$, $C_f$, $H_0$ such that
\begin{equation}\label{moseresti}
|v(t)|_{L^\infty(\Omega)}\leq C^*\max\left\{1, V, V_\Gamma, E_0\right\}\quad
\hbox{for a.e.~}t>0.
\end{equation}
\end{proposition}


We now finish the proof of Theorem \ref{main} by showing that $T_0$ introduced at the beginning
of this subsection is $+\infty$ if $R$ is sufficiently large.
Using (\ref{nequ2}), we obtain that
\begin{eqnarray}\label{me10}
|U(x,t)| &\le& C_{8}\left(1 + \int_0^t \expe^{\tau - t} \theta(x,\tau) \dd\tau\right)
\quad\mae\,,\\ \label{me11}
|U_t(x,t)| &\le& C_{9}\left(1 + \theta(x,t)
+ \int_0^t \expe^{\tau - t} \theta(x,\tau) \dd\tau\right)\quad\mae\,,
\end{eqnarray}
hence also (cf.~\eqref{nequ3})
\begin{equation}\label{me12}
|\chi_t(x,t)|\ \le\ C_{10}\left(1 + \theta(x,t)
+ \int_0^t\expe^{\tau - t}\theta(x,\tau)\dd\tau\right)\quad\mae
\end{equation}
By (\ref{esti1}), the function $U$ is in $L^\infty(0,\infty; L^2(\Omega))$ and the bound
does not depend on~$R$. Eq.~(\ref{nequ1}), with $\theta$ added to both the left and
the right hand side, thus satisfies the hypotheses of Proposition \ref{moser} for $N=3$ and $q=2$.
This enables us to conclude that $\theta(x,t)$ is uniformly bounded from above by a constant,
independently of $R$, so that $\theta$ never reaches the value $R$ if $R$ is sufficiently large,
which we wanted to prove. By (\ref{me10})--(\ref{me12}), also
$U$, $U_t$, and $\chi_t$ are uniformly bounded by a constant.

We proceed similarly to prove a uniform positive lower bound for $\theta$.
Set $R_0 := \sup\theta$, and in Eq.~(\ref{trequ1}) with $R>R_0$ put
$w = -\tilde w/\theta$, $\tilde w\in W^{1,2}(\Omega)$. For a new (nonnegative) variable $v(x,t) := \log R_0 - \log\theta(x,t)$
we obtain the equation
\begin{eqnarray}\label{me13}
&&\hspace{-15mm}\io v_t \tilde w(x)\dd x + \io \nabla v \cdot \nabla \tilde w(x) \dd x
+ \ipo h(x)\left(\frac{\theta_\Gamma}{\theta}-1\right) \tilde w(x) \dd s(x) \\ \nonumber
&=& \io\left(- \frac{U_t^2 +\chi_t^2}{\theta} -  \frac{|\nabla\theta|^2}{\theta^2}
+ U_t+ 2 \chi_t\right)  \tilde w(x)\dd x\,.
\end{eqnarray}
We now set
$$
\mathcal{H}[v] = \sign(v) \left(-\frac{U_t^2 +\chi_t^2}{\theta}-\frac{|\nabla\theta|^2}{\theta^2}
+ U_t+ 2 \chi_t\right)
$$
and check that the hypotheses of Proposition \ref{moser} are satisfied with
$f(v) = (\theta_\Gamma/R_0)(\expe^v - 1)$, $v_\Gamma = (R_0 -\theta_\Gamma)/R_0$,
$r \equiv 1$, and $v\mathcal{H}[v] \le 2C |v|$, where $C$ is a common upper bound for
$U_t$ and $\chi_t$. Hence, $v$ is bounded above by some $v^*$, which entails
$\theta \ge R_0 \expe^{-v^*}$.
This concludes the proof of Theorem \ref{main}.


\section{Long time behavior}\label{long}

We have the following statement.

\begin{proposition}\label{asp1}
Let the hypotheses of Theorem \ref{main} hold. Then we have
\begin{eqnarray}\label{ase1}
\int_0^\infty \left(\io \left(\theta_t^2 + U_t^2 + \chi_t^2 + |\nabla \theta|^2\right) \dd x
+ \ipo h(x)(\theta - \theta_\Gamma)^2\dd s(x)\right)\dd t &<& \infty,\qquad\\[2mm]\label{ase2}
\lim_{t\to\infty} \left(\io \left( U_t^2 + \chi_t^2 + |\nabla \theta|^2\right)(x,t) \dd x
+ \ipo h(x)(\theta - \theta_\Gamma)^2(x,t)\dd s(x)\right) &=& 0\,.
\end{eqnarray}
Furthermore, letting $t$ tend to $\infty$, the temperature $\theta$
converges strongly in $W^{1,2}(\Omega)$
to its equilibrium value $\theta_\Gamma$, and also both $\chi(x,t)$
and $U(x,t)$ converge strongly in $L^1(\Omega)$ (hence, strongly
in every $L^p(\Omega)$
for $p<\infty$) to their respective unique equilibrium values
$\chi_\infty$ and $U_\infty$ defined in Section \ref{mode}.
\end{proposition}

\bpf{Proof}\
The proof of the relations (\ref{ase1})--(\ref{ase2}) exactly follows
the argument of \cite[Proposition~5.1]{krsbottle1}.
For the convergence of the whole trajectory, it seems easier to prove it
directly without referring to the general theory
of dynamical systems with a unique equilibrium (see, e.g., \cite[Lemma~2.4]{chis}).

We eliminate from (\ref{sys1})--(\ref{sys3}) the term $U - \alpha (1-\chi)$ and obtain, using the notation (\ref{defm})--(\ref{equi4}),
the inclusion
\begin{equation}\label{long1}
A(x,t) + \frac{\alpha\lambda}{|\Omega|}
\io (\chi_\infty-\chi)\dd x' + \varrho_0 g (m + G_0\ell Z - x_3)
\in \partial I(\chi)\,,
\end{equation}
where we set
$$
A(x,t) := \nu U_t - \frac{\gamma}{\alpha} \chi_t + \Big(\frac{L}{\alpha\theta_c}
- \beta\Big) (\theta - \theta_\Gamma) + \frac{\beta}{|\Omega|}
\io (\theta - \theta_\Gamma)\dd x'\,.
$$
From (\ref{long1}) it follows that
\begin{equation}\label{long2}
\frac{\alpha\lambda}{|\Omega|}(\chi_\infty-\chi)
\io (\chi_\infty-\chi)\dd x' + \varrho_0 g(\chi_\infty-\chi)
(m + G_0\ell Z - x_3) \le (\chi_\infty-\chi) A(x,t)
\end{equation}
for a.e. $(x,t) \in \Omega_\infty$. Integrating (\ref{long2}) over
$\Omega$ and using (\ref{equi3a}) we obtain
\begin{equation}\label{long3}
\frac{\alpha\lambda}{|\Omega|}
\left|\io (\chi_\infty-\chi)\dd x'\right|^2
+ \varrho_0 g \io |\chi_\infty-\chi|\,
|m + G_0\ell Z - x_3'|\dd x' \le \io |\chi_\infty-\chi| |A(x',t)| \dd x'\,.
\end{equation}
By virtue of (\ref{ase2}), we have $\lim_{t\to \infty}
\io |A(x',t)| \dd x' = 0$. Hence, the right hand side of (\ref{long3})
converges to $0$ as $t \to \infty$. Since $\chi$ is a priori bounded,
the weighted $L^1$ convergence of $\chi_\infty - \chi$ with weight
$|m + G_0\ell Z - x_3|$ implies that $\chi \to \chi_\infty$ strongly
in $L^1(\Omega)$ as $t \to \infty$. The convergence $U \to U_\infty$
follows directly from (\ref{sys1}), (\ref{sys1eq}), and (\ref{ase2}).
\epf


\begin{thebibliography}{99}
{\small




\bibitem{bs}
  M.~Brokate, J.~Sprekels,
  Hysteresis and Phase Transitions,
  Appl. Math. Sci. 121,
  Springer,
  New York, 1996.

\bibitem{chis}
P.~Colli, D.~Hilhorst, F.~Issard-Roch, G.~Schimperna,
Long time convergence for a class of variational phase field models,
{\it ArXiv preprint no. 0801.2658v1\/} (2008).

\bibitem{fremond}
   M.~Fr\'emond,
    Non-smooth thermo-mechanics,
     Springer-Verlag Berlin, 2002.

\bibitem{fr1} M.~Fr\'emond, E.~Rocca, Well-posedness of a phase transition model
with the possibility of voids,
Math. Models Methods Appl. Sci.,
{\bf 16} no. 4 (2006),
559--586.

\bibitem{fr2} M.~Fr\'emond, E.~Rocca,
Solid liquid phase changes with different densities,
Q. Appl. Math., {\bf 66} (2008),
609--632.




\bibitem{hys}  P. Krej\v c\'\i, Hysteresis
operators -- a new approach to evolution differential
inequalities. {\it Comment. Math. Univ. Carolinae}, {\bf 33} no.~3
(1989), 525--536.



\bibitem{krsbottle1} P.~Krej\v c\'{\i}, E.~Rocca, J.~Sprekels,
A bottle in a freezer,
{\it ArXiv preprint no. 0904.4380v1\/} (2009).

\bibitem{krsplast} P.~Krej\v c\'{\i}, E.~Rocca, J.~Sprekels,
Phase changes in an elastoplastic container. In preparation.









\bibitem{visintin}
    A.~Visintin,
    Models of Phase Transitions,
     Progress in Nonlinear Differential Equations and their Applications 28,
     Birkh\"auser Boston, 1996.


\bibitem{ye} H.~Yildirim Erbil,
Surface Chemistry of Solid and Liquid Interfaces,
Blackwell Publishing, John Wiley \& Sons, 2006.




}
\end{thebibliography}
\end{document}